\documentclass[twocolumn]{autart}    
\usepackage{indentfirst} 
\usepackage{graphicx}     
\usepackage{amssymb}
\usepackage{CJK}
\usepackage{bm}
\usepackage{subfigure}
\usepackage{multirow}
\usepackage{booktabs}
\usepackage{pifont}
\usepackage{textcomp}
\usepackage{algorithm}
\usepackage{caption}
\usepackage{algpseudocode}
\usepackage[dvips]{epsfig}    

\begin{document}
\begin{frontmatter}

\title{Dual Privacy Guarantees for Distributed Nash Equilibrium Seeking in Aggregative Games\thanksref{footnoteinfo}} 

\thanks[footnoteinfo]{This work was supported in part by the NSFC under Grant 62473199, 62503228, 62533013, the Postdoctoral Fellowship Program of CPSF under Grant GZB20250958, and Jiangsu Funding Program for Excellent Postdoctoral Talent. The material in this paper was not presented at any conference. Corresponding author: Qian Ma.  }

\author[Paestum]{Qingtan Meng}\ead{qtmeng0609@163.com},    
\author[Paestum]{Qian Ma}\ead{qma@njust.edu.cn} 

\address[Paestum]{School of Automation, Nanjing University of Science and Technology, Nanjing, 210094, China.}  


\begin{abstract}                          
This paper investigates the privacy-preserving distributed Nash equilibrium seeking problem for aggregative games. A novel differential privacy mechanism is designed by incorporating stochastic event-triggering with stochastic quantization, which provides strong privacy protection by obfuscating the temporal patterns of information exchange among players and quantizing the transmitted information at triggering instants. Based on this mechanism, a differentially private distributed Nash equilibrium seeking algorithm with dual randomness is proposed. By embedding a decaying factor sequence into both the triggering condition and interaction terms among players, it is proved that the proposed algorithm can achieve rigorous $(0,\delta)$-differential privacy at each iteration while maintaining provable convergence. Crucially, this privacy guarantee is sustained over infinite iterations for a sufficiently large quantization interval and a sufficiently small trigger threshold tuning coefficient. Moreover, the synergy between event-triggered communication and quantization significantly enhances communication efficiency. Simulation results verify the validity of the proposed approach.

\end{abstract}

\end{frontmatter}

\section{Introduction}

In recent years, distributed Nash equilibrium (NE) seeking for noncooperative games has attracted increasing attention due to its wide applications in various fields  including but not limited to smart grids \cite{Tushar}, mobile sensor networks \cite{Stankovic}, and wireless communication  \cite{Fan}. As a powerful model of noncooperative games, aggregative games effectively capture competitive scenarios with aggregative interactions among multiple decision-makers. In such games, each player aims to optimize its own objective function, which depends on its own decision and the aggregate of all players' decisions, by exchanging information with its neighbors via local communication networks. Some pioneering works include \cite{Koshal,Ye,liang} for  distributed Nash equilibrium seeking for aggregative games. Consequently,  fruitful works have emerged on this issue; see \cite{Zhu,deng,de,Liul,Feng,Carnevale} and the reference therein.
\\
\indent\setlength{\parindent}{1em}
One typical feature of most existing distributed NE seeking algorithms is the explicit sharing of decision/estimation variables at each iteration, which may raise privacy concerns for players. Given this privacy risk, it is crucial to protect the sensitive information of players from being compromised while achieving NE seeking.  In \cite{Lu} and \cite{ZhuM}, privacy-preserving distributed Nash equilibrium seeking algorithms were developed based on secure multi-party computation. By injecting relevant noise to protect the privacy of players in  \cite{Gade}, fully distributed Nash equilibrium seeking for aggregative games was achieved. The authors in \cite{Shakarami} employed uncertain parameters to obfuscate the pseudo-gradient mapping, thereby addressing privacy-preserving distributed Nash equilibrium seeking in a continuous-time setting. In \cite{cai}, a privacy-preserving distributed algorithm was proposed via masked interactive information, ensuring that players' decisions converge to the Nash equilibrium. Neverthless, the aforementioned approaches only offer qualitative privacy guarantees, and cannot precisely quantify the protection strength for players' sensitive information. 
\\
\indent\setlength{\parindent}{1em}
Differential privacy has gained prominence by virtue of its ability to provide quantifiable and strong privacy protection against arbitrary postprocessing and auxiliary information \cite{Dwork}. By introducing uncorrelated random noise into the information exchanged among players, differentially private distributed NE seeking for aggregative games was first implemented in \cite{Ye2021}. On this basis, the authors in \cite{wang2022} further ensured NE convergence for stochastic aggregative games while achieving $(\epsilon,\delta)$-differential privacy.  In \cite{Zhang2025}, a differentially private distributed NE seeking algorithm was proposed based on gradient tracking and communication compression. In \cite{WangL2024}, distributed computation algorithms were developed for NE seeking  in linear-quadratic network games with provable differential privacy guarantees. By injecting Laplacian noise free from the requirement of uniformly bounded gradient information in \cite{Guo2024}, $\epsilon$-differential privacy was achieved  for the proposed distributed NE seeking strategies.  Although differentially private distributed NE seeking has made encouraging progress, all the aforementioned approaches suffer from a fundamental trade-off between privacy level and convergence accuracy, i.e., strong privacy protection and precise convergence cannot be guaranteed simultaneously.  To overcome this drawback, \cite{Wang2024} incorporated a weakening factor into the $\epsilon$-differentially private distributed NE seeking algorithm, ensuring convergence to NE despite the presence of persistent noise. Since then, several distributed NE seeking algorithms that can simultaneously achieve both rigorous differential privacy and guaranteed convergence have been developed in \cite{Chen2025,Wang20241,Odeyomi}.
\\
\indent\setlength{\parindent}{1em}
Recently, stochastic quantization effects have been recognized as a promising tool for enabling differential privacy, with the seminal work by \cite{Wang2022} demonstrating that a stochastic ternary quantization scheme can achieve $(0,\delta)$-differential privacy  in distributed optimization. Unlike conventional noise injection mechanisms, quantization-based differential privacy ensures privacy protection by compressing and randomizing the transmitted information, which also helps reduce communication overhead. Following this line of research, privacy-preserving distributed NE seeking algorithms  were proposed in \cite{Huo2024} and \cite{Zhao2025}, both achieving $(0,\delta)$-differential privacy through a stochastic uniform quantization scheme.  However, a key limitation of the existing works \cite{Wang2022,Huo2024,Zhao2025}  is that they only guarantee $(0,\delta)$-differential privacy at each iteration. As a result, $(0,1)$-differential privacy can be achieved over an infinite number of iterations. This degradation implies that the proposed algorithms in \cite{Wang2022,Huo2024,Zhao2025} would eventually output sensitive information directly, thus failing to provide any long-term privacy protection. 
\\
\indent\setlength{\parindent}{1em}
To address the issue of being unable to guarantee privacy over an infinite number of iterations,  we take the first step to propose a distributed NE seeking approach that can achieve strong privacy protection over infinite iterations without sacrificing convergence accuracy. This is accomplished by leveraging stochastic event-triggering with quantization effects, which provides a new perspective on differential privacy. The main contribution of the paper is as follows: 
\\
\indent\setlength{\parindent}{1em}
1) An innovative differential privacy mechanism with dual randomness is constructed.  More specifically, a stochastic event-triggering scheme is developed to obfuscate the temporal pattern of information exchange among players, making the transmission moments unpredictable to potential adversaries. Meanwhile, a stochastic quantization scheme is utilized to inject uncertainties into the triggered transmission values, preventing adversaries from inferring the original decision/estimate information from the observed data. 
\\
\indent\setlength{\parindent}{1em}
2) A novel differentially private distributed NE seeking algorithm for aggregative games is proposed.  By embedding a decaying factor sequence into both the triggering condition and the inter-player interaction terms, the proposed algorithm achieves rigorous $(0,\delta)$-differential privacy at each iteration without sacrificing provable convergence accuracy. Remarkably, while existing quantization enabled results only ensure $(0,1)$-differential privacy over infinite iterations \cite{Wang2022,Huo2024,Zhao2025}, our algorithm can maintain rigorous $(0,\delta)$-differential privacy over an infinite horizon, given a sufficiently large quantization interval and a sufficiently small trigger threshold tuning coefficient. 
\\
\indent\setlength{\parindent}{1em}
3) The synergistic effect of event-triggering and quantization schemes enables us to improve communication efficiency by reducing the transmission frequency and compressing the transmitted information. This advantage renders our approach more suitable for scenarios with limited communication resources than existing works on differential privacy distributed NE seeking \cite{Ye2021,wang2022,Zhang2025,WangL2024,Guo2024,Wang2024,Chen2025,Wang20241,Odeyomi,Wang2022,Huo2024,Zhao2025}.

\noindent\textbf{Notations.} $\mathbb{R}^i$ denotes Euclidean space with dimension $i$. 
For  matrix $P$, $\left\|P\right\|$ is its Euclidean norm and $P^T$ is its transpose.
For a real number $x$, $|x|$ denotes its absolute value.
The $N$-dimensional vectors $(1,\cdots,1)^T$ is represented by $\mathbf{1}_N$.
$\mathrm{diag}\{a_1,\cdots,a_n\}$ represents the diagonal matrix of elements $a_1,\cdots,a_n$. 
$\langle\cdot,\cdot\rangle$ denotes the inner product.

\section{Preliminaries and Problem Formulation}
\subsection{Aggregative Game}

In this paper, we consider the aggregative game problem among a networked of multiple players. All the definitions related to the game in the following are derived from \cite{Ye2021}.

\noindent\textbf{Definition 1 (A Normal Form Game).} A normal form game is defined as a triple $\Gamma=\{\mathcal{N},\Omega,\tilde{f}\}$, where $\mathcal{N}=\{1,2,\cdots,N\}$ is the set of players, $\Omega=\Omega_1\times\Omega_2\cdots\times\Omega_N$ with $\Omega_i\in\mathbb{R}$ being the decision set of player $i$, and  $\tilde{f}=(\tilde{f}_1,\tilde{f}_2,\cdots,\tilde{f}_N)$ with $\tilde{f}_i$ being the cost function of player $i$.

\noindent\textbf{Definition 2 (NE).} NE is an decision profile on which no player can reduce its cost by unilaterally changing its own action, i.e., an action profile $\mathbf{x}^*=(x_i^*,\mathbf{x}_{-i}^*)\in\Omega$ is a NE if for $i\in\mathcal{N}$, $\tilde{f}_i(x_i^*,\mathbf{x}_{-i}^*)\leq \tilde{f}_i(x_i,\mathbf{x}_{-i}^*)$, where $x_i\in\Omega_i$ and $\mathbf{x}_{-i}=[x_1,x_2,\cdots,x_{i-1},x_{i+1},\cdots,x_N]^T$.

Throughout this paper, we consider that there exists a function $f_i(x_i,\bar{x})$ such that $\tilde{f}_i(x_i,\mathbf{x}_{-i})=f_i(x_i,\bar{x})$, where $\bar{x}=\frac{1}{N}\sum_{i=1}^Nx_i$ denotes the aggregate average of all players's decisions. Recall that the decision variable of player $i$ is restricted in $\Omega_i$, the aggregate average $\bar{x}$ is restricted to the set $\bar{\Omega}=\frac{1}{N}(\Omega_1+\Omega_2+\cdots+\Omega_N)$ based on the Minkowski sum. With this notation, the aggregative game in which player $i$ acts can be indicated as 
\begin{eqnarray}
\min f_i(x_i,\bar{x}) \quad \mathrm{s.t.}\quad x_i\in\Omega_i \quad \mathrm{and}\quad \bar{x}\in\bar{\Omega}.\label{game}
\end{eqnarray}
Instead of relaying on a central authority, an undirected and connected graph $\mathcal{G}=(\mathcal{N},\mathcal{E})$ is used to described the information exchange among players, where $\mathcal{E}\subseteq \mathcal{N}\times \mathcal{N}$ is the edge set. For graph $\mathcal{G}$, $(j,i)\in\mathcal{E}$ represents that player $i$ can obtain information from player $j$. If $(j,i)\in\mathcal{E}$, then there exists $(i,j)\in\mathcal{E}$ since $\mathcal{G}$ is undirected. An undirected graph is connected if there is a path between each pair of players. The weight matrix of $\mathcal{G}$ is defined as $L=\{L_{ij}\}$ with $L_{ij}>0$ if $(j,i)\in\mathcal{E}$ and $L_{ij}=0$ otherwise. For player $i\in\mathcal{N}$, its neighbor set $\mathbb{N}_i$ is the collection of players $j$ with $L_{ij}>0$ and we define $L_{ii}=-\sum_{j\in\mathbb{N}_i}L_{ij}$. It is obvious that $\mathbf{1}^TL=\mathbf{0}^T$ and $L\mathbf{1}=\mathbf{0}$.

In this article, we  investigate the distributed NE seeking for aggregative game (\ref{game}). Namely, no player has direct access to the average decision $\bar{x}$. Instead, each player needs to construct a local estimate of the average value through local interactions with its neighbors. Setting $F(\mathbf{x},u)=[F_1(x_1,u),F_2(x_2,u),\cdots,F_N(x_N,u)]^T$ with $F_i(x_i,u)=\nabla_{x_i} f_i(x_i,u)$, the following assumptions are utilized to establish the main results of this article.

\noindent\textbf{Assumption 1.} For each player $i\in\mathcal{N}$, the decision set $\Omega_i$ is compact and convex, and the function $f_i(x_i,u)$ is convex in $x_i$ over the set $\Omega_i$, while $f_i(x_i,u)$ is continuous differentiable in $(x_i,u)$ over some open set containing the set $\Omega_i\times \bar{\Omega}$.

\noindent\textbf{Assumption 2.} The mapping $F_i(x_i,u)$ is Lipschitz continuous with respect to $u$, i.e., there exists a Lipschitz constant $l>0$ such that $|F_i(x_i,u_1)-F_i(x_i,u_2)|\leq l |u_1-u_2|$ for any fixed $x_i\in\Omega_i$ and all $u_1$, $u_2\in\bar{\Omega}$.
 
\noindent\textbf{Assumption 3.}  For $\mathbf{x}=[x_1,x_2,\cdots,x_N]^T$, the mapping $\phi(\mathbf{x})=F(\mathbf{x},\bar{x})$ is strictly monotone over $\Omega$, i.e.,  for all $\mathbf{x}\ne \mathbf{x}'$, the inequality $(\phi(\mathbf{x})-\phi(\mathbf{x}'))^T(\mathbf{x}-\mathbf{x}')>0$ always holds.

\noindent\textbf{Remark 1.} Assumptions 1 and 2 are fairly standard for the distributed NE seeking \cite{liang,Zhu}. The strictly monotone condition of Assumption 3 is weaker than the commonly strongly monotone assumption in \cite{Ye2021,Wang2022}. Under Assumptions 1 and 3, it follows from \cite{bianfen} that there is a unique NE and players' decisions are at the NE if and only if $\phi(\mathbf{x})=\textbf{0}_N$. Moreover, it follows from Assumption 1 that the set $\bar{\Omega}$ is convex since $\Omega_i$ is convex., i.e., the decision profile has an upper bound. On the basis, there exists a constant $\eta>0$ such that $|F_i(\cdot)|\leq \eta$.

\subsection{Differential Privacy}
Note that networked players need to interact with their neighbors, which may lead to the leakage of their sensitive information. In this paper, we propose to leverage quantization and  event-triggering effects to enable differential privacy in distributed NE seeking problem of aggregative game (\ref{game}). First, we define adjacency function sets for the problem.

 \noindent\textbf{Definition 3 (Adjacent).} Two function sets $\mathbf{F}=\{f_i\}_{i=1}^N$ and  $\mathbf{F^{'}}=\{f_i^{'}\}_{i=1}^N$ are said to be adjacent if the following conditions hold: 
 \\
\indent\setlength{\parindent}{1em}
 1) There exists an $i\in\mathcal{N}$ such that $f_i\ne f_i^{'}$ but $f_j=f_j^{'}$ for all $j\in\mathcal{N}$ and $j\ne i$;
 \\
\indent\setlength{\parindent}{1em}
 2) The different cost functions $f_i$ and $f_i^{'}$ have similar behaviors around $\mathbf{x}^*$, the NE of (\ref{game}) under the set $\mathbf{F}=\{f_i\}_{i=1}^N$. More specifically, there exists some $\alpha>0$ such that for all $\iota$ and $\iota^{'}$ in $B_{\alpha}(x_i^*)\triangleq\{v:v\in\mathbb{R},|v-x_i^*|<\alpha\}$, we have $\Pi_{\Omega_i}[\iota-\lambda\nabla_{\iota}f_i(\iota,\cdot)]-\iota=\Pi_{\Omega_i^{'}}[\iota^{'}-\lambda\nabla_{\iota^{'}}f_i^{'}(\iota^{'},\cdot)]-\iota^{'}$ for all $\lambda>0$, where $\Pi_{\Omega_i}[\cdot]$ denotes the Euclidean projection of a value onto the set $\Omega_i$.

In Definition 3, since the change of the cost function from $f_i$ to $f_i^{'}$ in the first condition can be arbitrary, certain restrictions must be imposed to ensure rigorous differential privacy in the distributed NE seeking procedure. Similar to \cite{Wang2024}, we adopt the second condition to guarantee rigorous differential privacy while maintaining provable convergence to an exact NE. It is worth pointing out that in the constraint-free case, this condition simplifies to requiring $\nabla_{\iota}f_i(\iota,\cdot)=\nabla_{\iota^{'}}f_i^{'}(\iota^{'},\cdot)$ for $\iota$ and $\iota^{'}$ in the neighborhood of the NE of (\ref{game}) under the cost function set $\mathbf{F}=\{f_i\}_{i=1}^N$.

Then, the definition of $(\epsilon,\delta)$-differential privacy is presented by following standard conventions \cite{Dwork}.

 \noindent\textbf{Definition 4 (Differential Privacy).}  Let $M$ be a distributed NE seeking algorithm. Given $\epsilon\geq 0$ and $\delta\geq 0$, the algorithm $M$ preserves $(\epsilon,\delta)$-differential privacy if for any pair of adjacent cost function sets $\mathbf{F}$, $\mathbf{F^{'}}$ and any $\mathcal{O}\subseteq Range(M)$, we always have
 \begin{eqnarray}\mathbb{P}(M(\mathbf{F})\in \mathcal{O})\leq e^{\epsilon}\mathbb{P}(M(\mathbf{F^{'}})\in\mathcal{O})+\delta,\nonumber
 \end{eqnarray} where $Range(M)$ is the domain of the observation under $M$. Note that $\mathcal{O}$ is the set of observations, i.e., the information transmitted in the network.

Definition 4 indicates that from the observed data sequence, it is almost impossible for attackers to distinguish these two function sets with a high probability. That is, it is difficult for attackers to obtain players' privacy. From Definition 4, a smaller $\epsilon\geq 0$ or $\delta\geq 0$ means higher privacy level.

\section{Differentially Private Distributed NE Seeking}

In this section, we propose a differentially private distributed NE seeking algorithm for aggregative game (\ref{game}). 

\subsection{Privacy-Preserving Mechanism}

To deal with the privacy issue caused by the information exchange  of the estimates of average value, we develop a privacy-preserving mechanism, which consists of two components: (i) a stochastic event-trigger for randomizing the timing of communication; (ii) a stochastic quantizer that randomizes the transmitted values. The former provides temporal privacy, while the latter provides numerical privacy.

For player $i\in\mathcal{N}$ and iteration $k$, a stochastic event-trigger is defined as 
\begin{eqnarray}
 \varsigma_i^k=\left\{\begin{array}[]{rcl}
  1,  \quad \mathrm{if}\quad \xi_i^k>\sigma e^{-c(\rho_i^k)^2/\gamma^k} \\
   0,\quad \mathrm{otherwise}   \qquad\qquad\quad\  \\
    \end{array}\right.,\label{trigger}
\end{eqnarray} 
where $\sigma> 1$ is a design parameter, $c>0$ is called a trigger threshold tuning coefficient, $\xi_i^k$ is a stationary ergodic random process taking values in $(a,1)$ with $a> 0$. For simplicity, we assume $\xi_i^k\sim \mathrm{Uniform}(a,1)$ in the following analysis. $\gamma^k>0$ is a decaying factor which is crucial for ensuring the convergence accuracy of the algorithm, $\rho_i^k$ is referred to as a trigger function which satisfies
\begin{eqnarray}
\rho_i^k=\tilde{y}_i^{\tau_i(k-1)}-y_i^k,\label{tf}
\end{eqnarray} 
where $y_i^k$ is the estimate of the average aggregate decision, $\tau_i(k-1)$ is a triggering time indicator that takes the form 
\begin{eqnarray}
\tau_i(k-1)=\left\{\begin{array}[]{rcl}
  k-1,  \qquad \mathrm{if}\quad \varsigma_i^{k-1}=1 \\
   \tau_i(k-2),\quad \mathrm{if}\quad \varsigma_i^{k-1}=0
    \end{array}\right.,\label{if}
\end{eqnarray}
where $k=1,2,\cdots$. In particular, we set $\tau_i(-1)=0$. At each iteration $k$, define  local storage value $\tilde{y}_i^k$ as
 \begin{eqnarray}
\tilde{y}_i^k=\left\{\begin{array}[]{rcl}
  Q(y_i^k),  \quad \mathrm{if}\quad \varsigma_i^k=1 \\
   \tilde{y}_i^{k-1},\quad  \mathrm{if}\quad \varsigma_i^k=0
    \end{array}\right.,\label{storage}
\end{eqnarray}
where $Q(\cdot)$ is a stochastic quantizer denoted by 
\begin{eqnarray}
\left\{\begin{array}[]{rcl}
  \mathbb{P}(Q(z+nd)=nd)=1-\frac{z}{d}\ \ \\
  \mathbb{P}(Q(z+nd)=(n+1)d)=\frac{z}{d}
    \end{array}\right.,\label{quantizer}
\end{eqnarray}
where  $n\in\mathbb{Z}$, $z\in[0,d]$, $d>0$ is the quantization interval. It should be pointed out that for any real number $b\in(nd,(n+1)d]$, it can be expressed as $b=z+nd$. By defining the quantization error as $\omega_b=Q(b)-b$,
it is derived that $\mathbb{E}[\omega_b]=\mathbb{E}[Q(z+nd)]-(z+nd)=(nd(1-\frac{z}{d})+(n+1)d\frac{z}{d})-(z+nd)=0$, implying the unbiasedness of (\ref{quantizer}). Then, the variance of (\ref{quantizer}) is bounded, i.e., $\mathbb{E}[\omega_b^2]=\mathbb{E}[(Q(z+nd))^2]-(z+nd)^2=(nd)^2(1-\frac{z}{d})+((n+1)d)^2\frac{z}{d}-(z+nd)^2\leq \frac{d^2}{4}$.

In the mechanism, suppose that if $\varsigma_i^k=1$, player $i$ will broadcast $Q(y_i^k)$ to player $j\in\mathbb{N}_i$ and update the local storage value $\tilde{y}_i^k=Q(y_i^k)$. Otherwise, no value will be sent and keep $\tilde{y}_i^k=\tilde{y}_i^{k-1}$. Under such circumstances, $\tilde{y}_i^{\tau_i(k-1)}$ in (\ref{tf}) represents the latest value broadcast by player $i$ at iteration $k$.


 \noindent\textbf{Remark 2.} The proposed stochastic mechanisms preserve the essential rules of their deterministic counterparts. For stochastic event-trigger (\ref{trigger}), taking $\xi_i^k\sim \mathrm{Uniform}(a,1)$ yields the triggering probability
 \begin{eqnarray}
\mathbb{P}(\varsigma_i^k=1)=\frac{1}{1-a}(1-\max\{a,\sigma e^{-c(\rho_i^k)^2/\gamma^k}\}).\nonumber
\end{eqnarray}
Recall that $\rho_i^k$ represents the difference between player $i$' current value and its last broadcast  value. Thus,  $\mathbb{P}(\varsigma_i^k=1)$ increases monotonically with the absolute trigger error $|\rho_i^k|$ for small values and reaches unity when $|\rho_i^k|$ is large enough. This aligns with the law of deterministic event-trigger. For stochastic quantizer (\ref{quantizer}), when the value $y_i^k$ lies in $(nd,(n+1)d]$, it is quantized to the nearest endpoint with higher probability, which is consistent with the behavior of a deterministic uniform quantizer. In fact, when $\xi_i^k$ is a strictly positive constant, the stochastic event-trigger degenerates to a deterministic one. The stochastic quantizer degenerates to a deterministic uniform quantizer when the rounding rule becomes deterministic (e.g., always rounding to the nearest endpoint). Therefore, (\ref{trigger}) and (\ref{quantizer}) can be viewed as generalized versions of their deterministic counterpart. 

\subsection{Algorithm Design}

\begin{algorithm}[htp]
\caption{Differentially Private Distributed Algorithms for Aggregative Games}
\label{algorithm1}
\raggedright  
\begin{algorithmic}
\State \textbf{Initialization:} Iteration stepsize $\lambda^k$, decaying factor $\gamma^k$, weighted matrix $L$, quantization interval $d$, trigger parameters $\sigma$, $c$, $a$. Initialize $y_i^0 = x_i^0$.

\State \textbf{Iterations:} For $k = 0, 1, \ldots$, each player $i \in \mathcal{N}$ does:

\State  a) \textbf{Stochastic event-triggering:}

\If{$k = 0$}

    \noindent Set $\varsigma_i^0 = 1$ and $\tau_i(0) = 0$
\Else
   
    \noindent (1) Compute the triggering error $\rho_i^k = \tilde{y}_i^{\tau_i(k-1)} - y_i^k$;
   
    \noindent (2) Generate a random number $\xi_i^k \sim \mathrm{Uniform}(a, 1)$;
    
    \noindent (3) Check the triggering condition $\xi_i^k > \sigma e^{-c(\rho_i^k)^2 / \gamma^k}$. If satisfied, set $\varsigma_i^k = 1$ and update $\tau_i(k) = k$; otherwise, set $\varsigma_i^k = 0$ and keep $\tau_i(k) = \tau_i(k-1)$.
\EndIf

\State b) \textbf{Stochastic quantization and communication:}

\State (1) If $\varsigma_i^k = 1$, compute $Q(y_i^k)$ and broadcast it to all $j \in \mathcal{N}_i$, and update $\tilde{y}_i^k = Q(y_i^k)$. Otherwise, send no message and keep $\tilde{y}_i^k = \tilde{y}_i^{k-1}$.
\State (2) If a new message $Q(y_j^k)$ is received from $j \in \mathcal{N}_i$, update $\tilde{y}_j^k = Q(y_j^k)$. If no message is received, keep $\tilde{y}_j^k = \tilde{y}_j^{k-1}$.

\State  c) \textbf{Decision update:}
\begin{equation}
x_i^{k+1} = \Pi_{\Omega_i}\left[x_i^k - \lambda^k F_i(x_i^k, y_i^k)\right], \label{x}
\end{equation}
\begin{equation}
y_i^{k+1} = y_i^k + \gamma^k \sum_{j \in \mathcal{N}_i} L_{ij}\left(\tilde{y}_j^k - \tilde{y}_i^k\right) + x_i^{k+1}- x_i^k. \label{y}
\end{equation}
\end{algorithmic}
\end{algorithm}

A distributed NE seeking algorithm under the privacy-preserving mechanism is presented in Algorithm 1. At iteration $k = 0$, we bypass the event-triggering condition and forcibly set $\varsigma_i^0 = 1$ for all players.  It guarantees that every player broadcasts its initial quantized estimate $Q(y_i^0)$ to all neighbors, thereby ensuring that the algorithm can run properly.

 \noindent\textbf{Remark 3.} Algorithm 1 can be directly extended to deal with the circumstances in which the players' decisions are multidimensional. Nevertheless, the primary contribution of this work lies in the development of a novel privacy-preserving mechanism via stochastic event-triggering and quantization, rather than in handling the complexity of  multidimensional decisions. Therefore, for analytical simplicity and to highlight the core contribution of the privacy-preserving mechanism, we assume that the players' decisions are one-dimensional. 

\section{Convergence Analysis}

In this section, we will provide the convergence analysis of Algorithm 1. 

 \noindent\textbf{Theorem 1.} Under Assumptions 1-3, if there exists some $T\geq 0$ such that $\gamma^k$ and $\lambda^k$ satisfying the conditions as $\sum_{k=T}^{\infty}\gamma^k=\infty$, $\sum_{k=T}^{\infty}\lambda^k=\infty$,  $\sum_{k=T}^{\infty}(\gamma^k)^2<\infty$, and $\sum_{k=T}^{\infty}\frac{(\lambda^k)^2}{\gamma^k}<\infty$, then Algorithm 1 converges to the NE of (\ref{game}) almost surely.

\noindent\textbf{Proof.}
By defining the average of the estimates and the average of the decisions as $\bar{y}^k=\frac{1}{N}\sum_{i=1}^Ny_i^k$ and $\bar{x}^k=\frac{1}{N}\sum_{i=1}^Nx_i^k$, we firstly prove $\bar{y}^k=\bar{x}^k$, i.e.,  $\sum_{i=1}^Ny_i^k=\sum_{i=1}^Nx_i^k$. This proof is provided by induction.

For $k=0$, it follows from $y_i^0=x_i^0$ that $\sum_{i=1}^Ny_i^0=\sum_{i=1}^Nx_i^0$.
Then, suppose that $\sum_{i=1}^Ny_i^k=\sum_{i=1}^Nx_i^k$ holds for some iteration $k>0$. In the following, we prove that it also holds for iteration $k+1$. 

Based on (\ref{x}) and (\ref{y}), it is directly deduced that 
\begin{eqnarray}
\sum_{i=1}^Ny_i^{k+1}=\gamma^k\sum_{i=1}^N\sum_{j\in\mathbb{N}_i}L_{ij}(\tilde{y}_j^k-\tilde{y}_i^k)+\sum_{i=1}^Nx_i^{k+1}.
\end{eqnarray}
 From the symmetric property of $L_{ij}$ in the undirected and connected graph, the first term of the right-hand side can be expressed as $\gamma^k\sum_{i=1}^N\sum_{j\in\mathbb{N}_i}L_{ij}(\tilde{y}_j^k-\tilde{y}_i^k)=\gamma^k\sum_{i=1}^N\sum_{j\in\mathbb{N}_i}L_{ij}\tilde{y}_j^k-\gamma^k\sum_{i=1}^N\sum_{i\in\mathbb{N}_j}L_{ji}\tilde{y}_i^k=0$. This means that $\sum_{i=1}^Ny_i^{k+1}=\sum_{i=1}^Nx_i^{k+1}$. Therefore, $\bar{y}^k=\bar{x}^k$ for all $k\geq 0$. Then, under Proposition 1 of \cite{Wang2024}, we will separately analyze $\sum_{i=1}^N(y_i^{k+1}-\bar{y}^{k+1})^2$ and $\sum_{i=1}^N(x_i^{k+1}-x_i^*)^2$.

\noindent\textbf{Part I: Analyze \(\bm{\sum_{i=1}^N(y_i^{k+1}-\bar{y}^{k+1})^2}\).}

Note that $\tilde{y}_i^k$ remains unchanged at non-triggering moments. According to (\ref{if}) and (\ref{storage}), it is not difficult to obtain 
\begin{eqnarray}
\tilde{y}_i^k=\tilde{y}_i^{\tau_i(k)}=Q(y_i^{\tau_i(k)}).
\end{eqnarray}
 Moreover, the equation $y_i^{\tau_i(k)}=\varsigma_i^ky_i^k+(1-\varsigma_i^k)y_i^{\tau_i(k-1)}$ holds. Letting $\varepsilon_i^k=Q(y_i^{\tau_i(k)})-y_i^{\tau_i(k)}$, it can be derived that
\begin{eqnarray}
\tilde{y}_i^k=\varepsilon_i^k-(1-\varsigma_i^k)\varepsilon_i^{k-1}+(1-\varsigma_i^k)\rho_i^k+y_i^k.\label{tildey2}
\end{eqnarray}
Based on the fact that $L_{ii}=-\sum_{j\in\mathbb{N}_i}L_{ij}$, (\ref{y}) can be rewritten as 
\begin{eqnarray}
y_i^{k+1} 
= y_i^k+ \gamma^k \sum_{j=1}^NL_{ij}(\varepsilon_j^k+y_j^k+(1-\varsigma_j^k)(\rho_j^k-\varepsilon_j^{k-1}))\nonumber\\
+x_i^{k+1}- x_i^k.\qquad\qquad\qquad\qquad\qquad\qquad\qquad
\end{eqnarray}
Defining $\mathbf{y}^k=[y_1^k,\cdots,y_N^k]^T$, $\mathbf{x}^k=[x_1^k,\cdots,x_N^k]^T$,  $\mathbf{\varepsilon}^k=[\varepsilon_1^k,\cdots,\varepsilon_N^k]^T$, $\mathbf{\rho}^k=[\rho_1^k,\cdots,\rho_N^k]^T$, one yields 
\begin{eqnarray}
\mathbf{y}^{k+1}=\mathbf{y}^k+\gamma^kL(\mathbf{\varepsilon}^k+\mathbf{y}^k+\varsigma^k(\mathbf{\rho}^k-\mathbf{\varepsilon}^{k-1}))\nonumber\\
+\mathbf{x}^{k+1}-\mathbf{x}^k,\qquad\qquad\qquad\qquad\quad
\end{eqnarray}
 where $\varsigma^k=\mathrm{diag}\{1-\varsigma_1^k,\cdots,1-\varsigma_N^k\}$. Based on the definition of $\bar{y}^{k+1}$ and $\mathbf{1}^TL=0$, we have
 \begin{eqnarray}
\bar{y}^{k+1}=\frac{1}{N}\sum_{i=1}^{N}y_i^{k+1}=\frac{\mathbf{1}^T}{N}\mathbf{y}^{k+1}\nonumber\\
=\frac{\mathbf{1}^T}{N}(\mathbf{y}^k+\mathbf{x}^{k+1}-\mathbf{x}^k),\ \ \label{baryk+1}
\end{eqnarray}
which implies 
\begin{eqnarray}
\mathbf{y}^{k+1}-\mathbf{1}\bar{y}^{k+1}=A^k\mathbf{y}^k+B^k(\mathbf{x}^{k+1}-\mathbf{x}^k)\qquad\nonumber\\
+\gamma^kL(\mathbf{\varepsilon}^k+\varsigma^k(\mathbf{\rho}^k-\mathbf{\varepsilon}^{k-1})),
\end{eqnarray} 
where $A^k=I+\gamma^k-\frac{\mathbf{1}\mathbf{1}^T}{N}$ and $B^k=I-\frac{\mathbf{1}\mathbf{1}^T}{N}$.
It is obvious that $A^k\mathbf{1}=0$. Therefore, the preceding equation further deduces
\begin{eqnarray}
\mathbf{y}^{k+1}-\mathbf{1}\bar{y}^{k+1}=A^k(\mathbf{y}^k-\mathbf{1}\bar{y}^{k})+B^k(\mathbf{x}^{k+1}-\mathbf{x}^k)\nonumber\\
+\gamma^kL(\mathbf{\varepsilon}^k+\varsigma^k(\mathbf{\rho}^k-\mathbf{\varepsilon}^{k-1})),
\end{eqnarray} 
which derives
\begin{eqnarray}
\left\|\mathbf{y}^{k+1}-\mathbf{1}\bar{y}^{k+1}\right\|^2\qquad\qquad\qquad\qquad\qquad\qquad\qquad\nonumber\\
=\left\|A^k(y^k-\mathbf{1}\bar{y}^{k})+B^k(\mathbf{x}^{k+1}-\mathbf{x}^k)+\gamma^kL\varsigma^k\mathbf{\rho}^k\right\|^2\ \nonumber\\
+2\langle A^k(y^k-\mathbf{1}\bar{y}^{k})+B^k(\mathbf{x}^{k+1}-\mathbf{x}^k)+\gamma^kL\varsigma^k\mathbf{\rho}^k,\nonumber\\
\gamma^kL(\mathbf{\varepsilon}^k-\varsigma^k\mathbf{\varepsilon}^{k-1})\rangle+\left\|\gamma^kL(\mathbf{\varepsilon}^k-\varsigma^k\mathbf{\varepsilon}^{k-1})\right\|^2\quad\; \nonumber\\
\leq \left\|A^k(\mathbf{y}^k-\mathbf{1}\bar{y}^{k})+B^k(\mathbf{x}^{k+1}-\mathbf{x}^k)+\gamma^kL\varsigma^k\mathbf{\rho}^k\right\|^2\ \nonumber\\
+2\langle A^k(\mathbf{y}^k-\mathbf{1}\bar{y}^{k})+B^k(\mathbf{x}^{k+1}-\mathbf{x}^k)+\gamma^kL\varsigma^k\mathbf{\rho}^k,\nonumber\\
\gamma^kL(\mathbf{\varepsilon}^k-\varsigma^k\mathbf{\varepsilon}^{k-1})\rangle\qquad\qquad\qquad\qquad\qquad\qquad\nonumber\\
+2(\gamma^k)^2\left\|L\right\|^2(\left\|\mathbf{\varepsilon}^k\right\|^2+\left\|\mathbf{\varepsilon}^{k-1}\right\|^2),\qquad\qquad\qquad\
\end{eqnarray} 
where $\left\|\varsigma^k\right\|\leq 1$. Reviewing the mean  and variance of the quantizer (\ref{quantizer}), it can be verified that $\mathbb{E}[\varepsilon_i^k]=0$ and $\mathbb{E}[(\varepsilon_i^k)^2]\leq \frac{d^2}{4}$ at any iteration $k$.

Let $\mathcal{F}^k=\{\mathbf{y}^0,\cdots,\mathbf{y}^k\}$. Using the fact that  $\left\|B^k\right\|=1$ and there exists a $T\geq 0$ such that $0<\left\|A^k\right\|\leq 1-\gamma^k|\varrho_2|$ holds for $k\geq T$ (see Lemma 5 in \cite{Wang2024}), where $\varrho_2$ is the second largest eigenvalue of $L$. Taking the conditional expectation yields
\begin{eqnarray}
\mathbb{E}[\left\|\mathbf{y}^{k+1}-\mathbf{1}\bar{y}^{k+1}\right\|^2|\mathcal{F}^k]\qquad\qquad\qquad\qquad\qquad\qquad\nonumber\\
\leq \left\|A^k(y^k-\mathbf{1}\bar{y}^{k})+B^k(\mathbf{x}^{k+1}-\mathbf{x}^k)+\gamma^kL\varsigma^k\mathbf{\rho}^k\right\|^2\ \ \; \nonumber\\
+(\gamma^k)^2\left\|L\right\|^2Nd^2\qquad\qquad\qquad\qquad\qquad\qquad\quad \nonumber\\
\leq[(1-\gamma^k|\varrho_2|)\left\|y^k-\mathbf{1}\bar{y}^{k}\right\|+\left\|\mathbf{x}^{k+1}-\mathbf{x}^k\right\|\qquad\quad\  \,\;\nonumber\\
+\gamma^k\left\|L\right\|\left\|\mathbf{\rho}^k\right\|]^2+(\gamma^k)^2\left\|L\right\|^2Nd^2\qquad\qquad\qquad\nonumber\\
\leq (1-\gamma^k|\varrho_2|)\left\|\mathbf{y}^k-\mathbf{1}\bar{y}^{k}\right\|^2+\frac{2}{\gamma^k|\varrho_2|}\left\|\mathbf{x}^{k+1}-\mathbf{x}^k\right\|^2\nonumber\\
+\frac{2}{\gamma^k|\varrho_2|}(\gamma^k)^2\left\|L\right\|^2\left\|\mathbf{\rho}^k\right\|^2+(\gamma^k)^2\left\|L\right\|^2Nd^2.\quad\ \,\label{exp}
\end{eqnarray}
According to (\ref{trigger}), when no trigger occurs, i.e., $\varsigma_i^k=0$, we can arrive at 
\begin{eqnarray}
(\rho_i^k)^2\leq \frac{\gamma^k}{c}(\mathrm{ln}\sigma-\mathrm{ln}\xi_i^k),
\end{eqnarray} 
which means that
\begin{eqnarray}
\left\|\mathbf{\rho}^k\right\|^2\leq\sum_{i=1}^N\frac{\gamma^k}{c}\mathrm{ln}(\sigma-\xi_i^k)\leq \frac{N\gamma^k}{c}(\mathrm{ln}\sigma-\mathrm{ln}a).
\end{eqnarray} 
Then, it follows from $\left\|\mathbf{y}^k-\mathbf{1}\bar{y}^k\right\|^2=\sum_{i=1}^N(y_i^k-\bar{y}^k)^2$ and $\left\|\mathbf{x}^{k+1}-\mathbf{x}^k\right\|^2=\sum_{i=1}^N(x_i^{k+1}-x_i^k)^2$ that
\begin{eqnarray}
\mathbb{E}\left[\sum_{i=1}^N(y_i^{k+1}-\bar{y}^{k+1})^2|\mathcal{F}^k\right]\qquad\qquad\qquad\qquad\qquad\quad\nonumber\\
\leq[(1-\gamma^k|\varrho_2|)\sum_{i=1}^N(y_i^k-\bar{y}^k)^2+(\gamma^k)^2\left\|L\right\|^2Nd^2\quad\nonumber\\
+\frac{2}{\gamma^k|\varrho_2|}\sum_{i=1}^N(x_i^{k+1}-x_i^k)^2\qquad\qquad\qquad\qquad\qquad\nonumber\\
+\frac{2N(\gamma^k)^2\left\|L\right\|^2}{|\varrho_2|c}\mathrm{ln}\frac{\sigma}{a}.\qquad\qquad\qquad\qquad\qquad\quad\ \ \,\label{exp1}
\end{eqnarray}
Under Assumption 2, it follows from (\ref{x}) that
\begin{eqnarray}
|x_i^{k+1}-x_i^k|\qquad\qquad\qquad\qquad\qquad\qquad\qquad\quad\nonumber\\
=|\Pi_{\Omega_i}[x_i^k-\lambda^kF_i(x_i^k,y_i^k)]-x_i^k|\qquad\qquad\qquad\nonumber\\
\leq|x_i^k-\lambda^kF_i(x_i^k,y_i^k)-x_i^k|\qquad\qquad\qquad\qquad\ \nonumber\\
\leq\lambda^k|F_i(x_i^k,\bar{y}^k)-F_i(x_i^k,y_i^k)|+\lambda^k|F_i(x_i^k,\bar{y}^k)|\ \nonumber\\
\leq \lambda^k\eta+\lambda^kl|y_i^k-\bar{y}^k|.\qquad\qquad\qquad\qquad\qquad\ \
\end{eqnarray}
Therefore, we can obtain 
\begin{eqnarray}
(x_i^{k+1}-x_i^k)^2\leq 2\eta^2(\lambda^k)^2+2l^2(\lambda^k)^2(y_i^k-\bar{y}^k)^2.
\end{eqnarray}
Then, one has
\begin{eqnarray}
\mathbb{E}\left[\sum_{i=1}^N(y_i^{k+1}-\bar{y}^{k+1})^2|\mathcal{F}^k\right]\leq\frac{4N\eta^2(\lambda^k)^2}{\gamma^k|\varrho_2|}\qquad\nonumber\\
+\left(1-\gamma^k|\varrho_2|+\frac{4l^2(\lambda^k)^2}{\gamma^k|\varrho_2|}\right)\sum_{i=1}^N(y_i^k-\bar{y}^k)^2\nonumber\\
+(\gamma^k)^2\left\|L\right\|^2Nd^2+\frac{2N(\gamma^k)^2\left\|L\right\|^2}{|\varrho_2|c}\mathrm{ln}\frac{\sigma}{a}.\quad\label{part1}
\end{eqnarray}

\noindent\textbf{Part II: Analyze \(\bm{\sum_{i=1}^N(x_i^{k+1}-x_i^*)^2}\)}. 

At the NE $\mathbf{x}^*=[x_1^*,\cdots,x_N^*]^T$, the equation $x_i^*=\Pi_{\Omega_i}[x_i^*-\lambda^kF_i(x_i^*,\bar{x}^*)]$ always holds with $\bar{x}^*=\frac{1}{N}\sum_{i=1}^Nx_i^*$. Hence, it is derived that
\begin{eqnarray}
(x_i^{k+1}-x_i^*)^2\qquad\qquad\qquad\qquad\qquad\qquad\qquad\qquad\qquad\nonumber\\
=(\Pi_{\Omega_i}[x_i^k-\lambda^kF_i(x_i^k,y_i^k)]-\Pi_{\Omega_i}[x_i^*-\lambda^kF_i(x_i^*,\bar{x}^*)])^2\nonumber\\
\leq(x_i^k-\lambda^kF_i(x_i^k,y_i^k)-x_i^*+\lambda^kF_i(x_i^*,\bar{x}^*))^2\qquad\qquad\ \,\nonumber\\
=(x_i^k-x_i^*)^2\underbrace{-2\lambda^k(x_i^k-x_i^*)(F_i(x_i^k,y_i^k)-F_i(x_i^*,\bar{x}^*))}_{\textcircled{1}}\ \nonumber\\
\underbrace{+(\lambda^k)^2(F_i(x_i^k,y_i^k)-F_i(x_i^*,\bar{x}^*))^2}_{\textcircled{2}}.\label{x_i}\qquad\qquad\qquad\qquad
\end{eqnarray}
Under Assumption 2, it follows from $\bar{y}^k=\bar{x}^k$  that 
\begin{eqnarray}
\textcircled{1}
=-2\lambda^k(x_i^k-x_i^*)(F_i(x_i^k,y_i^k)-F_i(x_i^k,\bar{y}^k))\qquad\quad\nonumber\\
-2\lambda^k(x_i^k-x_i^*)(F_i(x_i^k,\bar{y}^k)-F_i(x_i^*,\bar{x}^*))\qquad\quad\nonumber\\
\leq \frac{(\lambda^k)^2}{\gamma^k}(x_i^k-x_i^*)^2+\gamma^k|F_i(x_i^k,y_i^k)-F_i(x_i^k,\bar{y}^k)|^2\nonumber\\
-2\lambda^k(x_i^k-x_i^*)(F_i(x_i^k,\bar{x}^k)-F_i(x_i^*,\bar{x}^*))\qquad\quad\nonumber\\
\leq\frac{(\lambda^k)^2}{\gamma^k}(x_i^k-x_i^*)^2+\gamma^kl^2(y_i^k-\bar{y}^k)^2\qquad\qquad\quad\ \,\nonumber\\
-2\lambda^k(x_i^k-x_i^*)(F_i(x_i^k,\bar{x}^k)-F_i(x_i^*,\bar{x}^*)).\qquad\quad\label{x_i1}
\end{eqnarray}
Then, we arrive at
\begin{eqnarray}
\textcircled{2}
\leq 2(\lambda^k)^2(|F_i(x_i^k,y_i^k)|^2+|F_i(x_i^*,\bar{x}^*)|^2)\nonumber\\
\leq 4\eta^2(\lambda^k)^2.\qquad\qquad\qquad\quad\ \qquad\qquad\label{x_i2}
\end{eqnarray}
Substituting (\ref{x_i1}) and (\ref{x_i2}) into (\ref{x_i}) yields
\begin{eqnarray}
\sum_{i=1}^N(x_i^{k+1}-x_i^*)^2\qquad\qquad\qquad\qquad\qquad\qquad\qquad\quad\nonumber\\
\leq (1+\frac{(\lambda^k)^2}{\gamma^k})\sum_{i=1}^N(x_i^k-x_i^*)^2+\gamma^kl^2\sum_{i=1}^N(y_i^k-\bar{y}^k)^2\nonumber\\
+4N\eta^2(\lambda^k)^2-2\lambda^k(\phi(\mathbf{x}^k)-\phi(\mathbf{x}^*))^T(\mathbf{x}^k-\mathbf{x}^*).\label{x_i3}
\end{eqnarray}
Recall from Assumption 3 that $(\phi(\mathbf{x}^k)-\phi(\mathbf{x}^*))^T(\mathbf{x}^k-\mathbf{x}^*)$  is always positive for all $\mathbf{x}^k\ne \mathbf{x}^*$. On the basis, by defining $\mathbf{z}^k=[\sum_{i=1}^N(x_i^k-x_i^*)^2,\sum_{i=1}^N(y_i^k-\bar{y}^k)^2]^T$, it is straightforward from Part I and Part II that
\begin{eqnarray}
\mathbb{E}[{\mathbf{z}}^{k+1}|\mathcal{F}^k]\qquad\qquad\qquad\qquad\qquad\qquad\qquad\qquad\qquad\nonumber\\
\leq 
    \left[
      \begin{array}{cc}
          1 & \gamma^kl^2  \\
          0 & 1-\gamma^k|\varrho_2|
     \end{array}
    \right]\mathbf{z}^k+\left[
      \begin{array}{cc}
          \frac{(\lambda^k)^2}{\gamma^k} & 0  \\
         0 & \frac{4l^2(\lambda^k)^2}{\gamma^k|\varrho_2|}
      \end{array}
    \right]\mathbf{z}^k\qquad\quad\ \nonumber\\
    +\left[
      \begin{array}{cc}
         4N\eta^2(\lambda^k)^2 \\
          (\gamma^k)^2\left\|L\right\|^2Nd^2+\frac{2\left\|L\right\|^2\mathrm{ln}\frac{\sigma}{a}(\gamma^k)^2}{c|\varrho_2|}+\frac{4N\eta^2(\lambda^k)^2}{\gamma^k|\varrho_2|} 
     \end{array}
    \right]\nonumber\\
    -2\lambda^k\left[
      \begin{array}{cc}
          (\phi(\mathbf{x}^k)-\phi(\mathbf{x}^*))^T(\mathbf{x}^k-\mathbf{x}^*)  \\
          0 
      \end{array}
    \right]\qquad\qquad\quad\ \ \nonumber\\
    \leq \left(\left[
      \begin{array}{cc}
          1 & \gamma^kl^2  \\
          0 & 1-\gamma^k|\varrho_2|
      \end{array}
    \right]+a^k\mathbf{1}\mathbf{1}^T\right)\mathbf{z}^k+b^k\mathbf{1}\qquad\qquad\quad \nonumber\\
    -\varkappa^k\left[
      \begin{array}{cc}
          (\phi(\mathbf{x}^k)-\phi(\mathbf{x}^*))^T(\mathbf{x}^k-\mathbf{x}^*)  \\
          0 
      \end{array}
    \right],\qquad\qquad\quad\ \
    \end{eqnarray}
where $\varkappa^k=2\lambda^k$, $a^k=\max\{\frac{(\lambda^k)^2}{\gamma^k},\frac{4l^2(\lambda^k)^2}{\gamma^k|\varrho_2|}\}$, and $b^k=\max\{4N\eta^2(\lambda^k)^2,(\gamma^k)^2\left\|L\right\|^2Nd^2+\frac{2N\left\|L\right\|^2\mathrm{ln}\frac{\sigma}{a}(\gamma^k)^2}{c|\varrho_2|}+\frac{4N\eta^2(\lambda^k)^2}{\gamma^k|\varrho_2|}\}$.
In light of the conditions that $\sum_{k=T}^{\infty}(\gamma^k)^2<\infty$ and $\sum_{k=T}^{\infty}\frac{(\lambda^k)^2}{\gamma^k}<\infty$,  it holds that $a^k$ and $b^k$ are summable.
Therefore, one gets that $\sum_{i=1}^N(x_i^k-x_i^*)^2$ and $\sum_{i=1}^N(y_i^k-\bar{y}^k)^2$ satisfy the conditions of Proposition 1 in \cite{Wang2024}. Hence, it can be deduced that $\lim_{k\to\infty}|x_i^k-x_i^*|=0$ and $\lim_{k\to\infty}|y_i^k-\bar{y}^k|=0$ almost surely for all $i\in\mathcal{N}$. The desired result follows.

 \noindent\textbf{Remark 4.} Theorem 1 ensures that Algorithm 1 achieves exact convergence to the NE. In contrast to existing algorithms that can only converge to a neighborhood of the NE, our design incorporates two key features. First, each player $i$ employs the quantized estimate $\tilde{y}_i^k$ (rather than its raw estimate $y_i^k$) in its interaction terms $L_{ij}(\tilde{y}_j^k-\tilde{y}_i^k)$. Second,  a time-diminishing sequence $\{\gamma^k\}$ is introduced in both event-triggering condition (\ref{trigger}) and the interaction terms of (\ref{y}). These designs effectively eliminate the influence of quantization errors and suppress the impact of uncertainties arising from the event-trigger on the average aggregate  estimation. 

 \noindent\textbf{Remark 5.} With the help of Lemma 4 in \cite{wangcaijian}, it follows from (\ref{part1}) that $\sum_{i=1}^N(y_i^{k+1}-\bar{y}^{k+1})^2$ decays to zero at a rate of $\mathcal{O}((\frac{\lambda^k}{\gamma^k})^2)$. Further, if the mapping $\phi(\mathbf{x})$ degenerates into strongly monotone, i.e., $(\phi(\mathbf{x})-\phi(\mathbf{x}'))^T(\mathbf{x}-\mathbf{x}')>\zeta\left\|\mathbf{x}-\mathbf{x}'\right\|^2$ for some $\zeta>0$, then it follows from (\ref{x_i3}) that $\sum_{i=1}^N(x_i^{k+1}-x_i^*)^2$ decays to zero at a rate of $\mathcal{O}(\frac{\lambda^k}{\gamma^k})$.

\section{Privacy Analysis}

Note that the output of Algorithm 1 is the estimate of the average aggregate decision, i.e., $\mathbf{y}^k=[y_1^k,\cdots,y_N^k]^T$. With the help of the idea of differential privacy design in \cite{Huang}, we define the sensitivity of the proposed algorithm  as follows.

 \noindent\textbf{Definition 5 (Sensitivity).}  At each iteration $k$, for any adjacent function sets $\mathbf{F}$ and $\mathbf{F^{'}}$, the sensitivity of Algorithm 1  is defined as $\Delta^k=\sup_{\mathrm{Adj}(\mathbf{F},\mathbf{F}^{'})}\left\|\mathbf{y}^k-\mathbf{y}^{'k}\right\|$.

It should be noted that Definition 5 is provided under any identical initial values (i.e., $x_i^0=x_i^{'0}=y_i^0=y_i^{'0}$ for $i\in\mathcal{N}$) of the prescribed adjacent cost function sets. Given the above definition, we have the following result.

 \noindent\textbf{Lemma 1.} At each iteration $k$, there exists a constant $\tilde{C}>0$ such that the sensitivity of Algorithm 1 satisfies 
 \begin{eqnarray}
 \Delta^k\leq \tilde{C}\frac{(\lambda^k)^2}{\gamma^k}.
 \end{eqnarray}

\noindent\textbf{Proof.}
Recall in Definition 3, we consider two adjacent function sets $\mathbf{F}=\{f_i\}_{i=1}^N$ and $\mathbf{F^{'}}=\{f_i^{'}\}_{i=1}^N$ where $f_i\ne f_i^{'}$ and $f_j=f_j^{'}$ for any $j\in\mathcal{N}$ and $j\ne i$.
It can be deduced that if the gradient $F_i=F_i^{'}$, the two sequences generated by the proposed algorithm under $\mathbf{F}$ and $\mathbf{F^{'}}$ would be the same by enforcing the event-trigger  (including trigger condition, trigger time, etc.) and quantizer (including quantization interval, etc.) to be the same, indicating that the proposed algorithm is of complete privacy. Hence, we only consider the case that $F_i\ne F_i^{'}$. Since the initial conditions, cost functions, and observations of $\mathbf{F}$ and $\mathbf{F^{'}}$ are identical for all $j\in\mathcal{N}$ and $j\ne i$, this implies that $y_j^k=y_j^{'k}$, i.e., $\left\|\mathbf{y}^k-\mathbf{y}^{'k}\right\|$ is always equal to $|y_i^k-y_i^{'k}|$.

Analogous to the sensitivity calculation in the noise-free case of standard differential privacy \cite{Wang2024}, we calculate the sensitivity without triggering and quantization, as the sensitivity is independent of the frequency of stochastic event-triggering and the quantization interval of stochastic quantizer for differential privacy. By (\ref{y}) of Algorithm 1, we have
\begin{eqnarray}
y_i^{k+1}-y_i^{'k+1}=(1-|L_{ii}|\gamma^k)(y_i^k-y_i^{'k})+x_i^{k+1}-x_i^k\nonumber\\
-(x_i^{'k+1}-x_i^{'k}).\label{yii}\qquad\qquad\qquad\qquad
\end{eqnarray}
Substituting (\ref{x}) into (\ref{yii}) yields
\begin{eqnarray}
|y_i^{k+1}-y_i^{'k+1}|\leq(1-|L_{ii}|\gamma^k)(y_i^k-y_i^{'k})\qquad\qquad\nonumber\\
+\lambda^k|F_i(x_i^k,y_i^k)-F_i^{'}(x_i^{'k},y_i^{'k})|.\label{yii2}
\end{eqnarray}
Note that when the conditions in the statement of Theorem 1 are satisfied, Algorithm 1 guarantees the convergence under both $\mathbf{F}$  and $\mathbf{F^{'}}$ to the respective NE, which are the same under the second requirement in Definition 3. This implies that $|F_i(x_i^k,y_i^k)-F_i^{'}(x_i^{'k},y_i^{'k})|=0$ will hold when $k$ is sufficiently large (for the iterates under  both $\mathbf{F}$  and $\mathbf{F^{'}}$ to enter the neighborhood $B_{\alpha}$ in Definition 3, upon which the evolution under $\mathbf{F}$  and $\mathbf{F^{'}}$ will be identical). Furthermore, the ensured convergence also means that $F_i(x_i^k,y_i^k)$ and $F_i^{'}(x_i^{'k},y_i^{'k})$ are always bounded. Therefore, there always exists some constant $C>0$ such that the following inequality  holds for all $k>0$ under the conditions of Theorem 1:
\begin{eqnarray}
|F_i(x_i^k,y_i^k)-F_i^{'}(x_i^{'k},y_i^{'k})|\leq C\lambda^k.
\end{eqnarray} 
 Then, (\ref{yii2}) can be rewritten as
 \begin{eqnarray}
|y_i^{k+1}-y_i^{'k+1}|\leq(1-|L_{ii}|\gamma^k)(y_i^k-y_i^{'k})+C(\lambda^k)^2.
\end{eqnarray}
 In light of Lemma 4 of \cite{wangcaijian}, it is straightforward that 
 \begin{eqnarray}
 |y_i^k-y_i^{'k}|\leq\tilde{C}\frac{(\lambda^k)^2}{\gamma^k},
 \end{eqnarray}
  where $\tilde{C}$ is some constant. Consequently, it can be deduced  $\Delta^k\leq \tilde{C}\frac{(\lambda^k)^2}{\gamma^k}$.

Now we show that based on the proposed privacy protection mechanism, our distributed NE seeking algorithm can achieve differential privacy while maintaining provable convergence.

 \noindent\textbf{Theorem 2.} Under Assumptions 1-3, Algorithm 1  is $(0,\delta^k)$-differentially private at each iteration $k$ if the design parameters satisfy
\begin{eqnarray}
 \left(\frac{\sigma}{1-a}\sqrt{\frac{2c}{e\gamma^k}}+\frac{1}{d}\right)\tilde{C}\frac{(\lambda^k)^2}{\gamma^k}=\delta^k\in(0,1).
\end{eqnarray}
 Moreover, when the sequence $\{\frac{(\lambda^k)^2}{(\gamma^k)^{3/2}}\}$ is summable, Algorithm 1 can achieve $(0,\sum_{k=0}^T\delta^k)$-differential privacy over $T$ iterations as $T\to\infty$ for a sufficiently large quantization interval and a sufficiently small trigger threshold tuning coefficient.

\noindent\textbf{Proof.}
Based on Definition 4, $M(\mathbf{F}^k)$ can be expressed as $M(\mathbf{F}^k)=[M_1^k,M_2^k,\cdots,M_N^k]^T$ at iteration $k$ with
\begin{eqnarray}
M_i^k=\left\{\begin{array}[]{rcl}
  Q(y_i^k),  \quad \mathrm{with\ \,prob.}\ \,\mathbb{P}(\varsigma_i^k=1|y_i^k) \\
  \emptyset,\quad \mathrm{with\ \,prob.}\ \,1-\mathbb{P}(\varsigma_i^k=1|y_i^k) 
      \end{array}\right., \  i\in\mathcal{N}.\nonumber
\end{eqnarray}
For player $i$ and iteration $k$, the proposed privacy-preserving mechanism implies that if the triggering condition (\ref{trigger}) is satisfied, the player will transmit the message $Q(y_i^k)$ to its neighbors. Otherwise, no message is sent, which is also observable to an eavesdropper. Therefore, for algorithms $M(\mathbf{F}^k)$ and $M(\mathbf{F^{'k}})$, the following conditions must be satisfied simultaneously to ensure that the observation sequences are identical: 
\\
\indent\setlength{\parindent}{1em}
1) The triggering decision $\varsigma_i^k$ and $\varsigma_i^{'k}$ should be identical. That is, $\varsigma_i^k=\varsigma_i^{'k}=1$ or $\varsigma_i^k=\varsigma_i^{'k}=0$; 
\\
\indent\setlength{\parindent}{1em}
2) If $\varsigma_i^k=\varsigma_i^{'k}=1$, then the transmitted messages should be identical, that is, $Q(y_i^k)=Q(y_i^{'k})$.

Next, we will proceed with the proof in two different cases.

\noindent\textbf{Case I:  \(\bm{\mathbf{\varsigma_i^k=\varsigma_i^{'k}=1}}\).} In this case, suppose $Q(y_i^k)=Q(y_i^{'k})=q$, which derives that 
\begin{eqnarray}
\mathbb{P}(M_i^k\in\mathcal{O}^k|y_i^k)=\mathbb{P}(\varsigma_i^k=1)\mathbb{P}(Q(y_i^k)=q|y_i^k),\ \; \nonumber\\
\mathbb{P}(M_i^{'k}\in\mathcal{O}^k|y_i^{'k})=\mathbb{P}(\varsigma_i^{'k}=1)\mathbb{P}(Q(y_i^{'k})=q|y_i^{'k}),
\end{eqnarray}
where $\mathcal{O}^k$ is the observation at iteration $k$.
Thus, one can compute that
\begin{eqnarray}
|\mathbb{P}(M_i^k\in\mathcal{O}^k|y_i^k)-\mathbb{P}(M_i^{'k}\in\mathcal{O}^k|y_i^{'k})|\qquad\qquad\qquad\quad\nonumber\\
\leq |\mathbb{P}(\varsigma_i^k=1)||\mathbb{P}(Q(y_i^k)=q|y_i^k)-\mathbb{P}(Q(y_i^{'k})=q|y_i^{'k})|\nonumber\\
+|\mathbb{P}(\varsigma_i^k=1)-\mathbb{P}(\varsigma_i^{'k}=1)||\mathbb{P}(Q(y_i^{'k})=q|y_i^{'k})|\qquad\nonumber\\
\leq |\mathbb{P}(Q(y_i^k)=q|y_i^k)-\mathbb{P}(Q(y_i^{'k})=q|y_i^{'k})|\qquad\qquad\ \;\nonumber\\
+|\mathbb{P}(\varsigma_i^k=1)-\mathbb{P}(\varsigma_i^{'k}=1)|.\qquad\qquad\qquad\qquad\qquad\,\label{p-p}
\end{eqnarray}
Noting that $\mathbb{P}(\varsigma_i^k=1)=\frac{1}{1-a}(1-\max\{a,\sigma e^{-c(\rho_i^k)^2/\gamma^k}\})$, it is not difficult to obtain 
\begin{eqnarray}
|\mathbb{P}(\varsigma_i^k=1)-\mathbb{P}(\varsigma_i^{'k}=1)|\qquad\qquad\qquad\nonumber\\
\leq \frac{\sigma}{1-a}\left|e^{-c(\rho_i^k)^2/\gamma^k}-e^{-c(\rho_i^{'k})^2/\gamma^k}\right|.
\end{eqnarray} 
It is straightforward from the mean value theorem that 
\begin{eqnarray}
\left|e^{-c(\rho_i^k)^2/\gamma^k}-e^{-c(\rho_i^{'k})^2/\gamma^k}\right|\leq\sqrt{\frac{2c}{e\gamma^k}}|\rho_i^k-\rho_i^{'k}|\qquad\nonumber\\
\leq\sqrt{\frac{2c}{e\gamma^k}}|\tilde{y}_i^{\tau_i(k-1)}-y_i^k-\tilde{y}_i^{'\tau_i(k-1)}+y_i^{'k}|\qquad\quad\nonumber\\
\leq\sqrt{\frac{2c}{e\gamma^k}}|Q(y_i^{\tau_i(k-1)})-y_i^k-Q(y_i^{\tau_i('k-1)})+y_i^{'k}|\nonumber\\
\leq\sqrt{\frac{2c}{e\gamma^k}}|y_i^k-y_i^{'k}|,\qquad\qquad\qquad\qquad\qquad\qquad\,
\end{eqnarray}
which leads to
\begin{eqnarray}
|\mathbb{P}(\varsigma_i^k=1)-\mathbb{P}(\varsigma_i^{'k}=1)|\leq \frac{\sigma}{1-a}\sqrt{\frac{2c}{e\gamma^k}}|y_i^k-y_i^{'k}|.\label{p-p1}
\end{eqnarray}
Then, without loss of generality, assume $y_i^k\in(nd,(n+1)d]$, i.e., $q=nd$ or $q=(n+1)d$. It follows from $Q(y_i^k)=Q(y_i^{'k})=q$ that   the value of $y_i^{'k}$ can be classified into three different scenarios:
\begin{enumerate}
    \item [$\bullet$] $y_i^{'k}\in (nd,(n+1)d]$;
    \item [$\bullet$] $y_i^{'k}\in ((n+1)d,(n+2)d]$;
    \item [$\bullet$] $y_i^{'k}\in ((n-1)d,nd]$.
    \end{enumerate}

For the first scenario, it can be obtained from (\ref{quantizer}) that 
\begin{eqnarray}
|\mathbb{P}(Q(y_i^k)=q|y_i^k)-\mathbb{P}(Q(y_i^{'k})=q|y_i^{'k})|\qquad\qquad\ \qquad \nonumber\\
=|\mathbb{P}(Q(y_i^k)=nd|y_i^k)-\mathbb{P}(Q(y_i^{'k})=nd|y_i^{'k})|\quad\ \;\qquad\nonumber\\
=\left|1-\frac{y_i^k-nd}{d}-(1-\frac{y_i^{'k}-nd}{d})\right|\leq \frac{|y_i^k-y_i^{'k}|}{d},\qquad\nonumber
\end{eqnarray}
or
\begin{eqnarray}
|\mathbb{P}(Q(y_i^k)=q|y_i^k)-\mathbb{P}(Q(y_i^{'k})=q|y_i^{'k})|\qquad\qquad\qquad\qquad \nonumber\\
=|\mathbb{P}(Q(y_i^k)=(n+1)d|y_i^k)-\mathbb{P}(Q(y_i^{'k})=(n+1)d|y_i^{'k})|\nonumber\\
=\left|\frac{y_i^k-nd}{d}-\frac{y_i^{'k}-nd}{d}\right|\leq \frac{|y_i^k-y_i^{'k}|}{d}.\qquad\qquad\quad\qquad\ \,\nonumber
\end{eqnarray}
Then, for the second scenario, it is straightforward that
\begin{eqnarray}
|\mathbb{P}(Q(y_i^k)=q|y_i^k)-\mathbb{P}(Q(y_i^{'k})=q|y_i^{'k})|\qquad\qquad\qquad\qquad \nonumber\\
=|\mathbb{P}(Q(y_i^k)=(n+1)d|y_i^k)-\mathbb{P}(Q(y_i^{'k})=(n+1)d|y_i^{'k})| \nonumber\\
=\left|\frac{y_i^k-nd}{d}-(1-\frac{y_i^{'k}-nd}{d})\right|\qquad\qquad\qquad\qquad\qquad\ \;\, \nonumber\\
\leq \frac{|y_i^k+y_i^{'k}-2(n+1)d|}{d}\leq \frac{|y_i^k-y_i^{'k}|}{d}.\qquad\qquad\qquad\quad\,  \nonumber
\end{eqnarray}
For the third scenario, we have
\begin{eqnarray}
|\mathbb{P}(Q(y_i^k)=q|y_i^k)-\mathbb{P}(Q(y_i^{'k})=q|y_i^{'k})|\qquad\qquad\ \qquad \nonumber\\
=|\mathbb{P}(Q(y_i^k)=nd|y_i^k)-\mathbb{P}(Q(y_i^{'k})=nd|y_i^{'k})|\quad\ \;\qquad\nonumber\\
=\left|1-\frac{y_i^k-nd}{d}-(1-\frac{y_i^{'k}-(n-1)d}{d})\right|\qquad\qquad\quad\,\nonumber\\
\leq \frac{|y_i^k+y_i^{'k}-2nd|}{d}\leq \frac{|y_i^k-y_i^{'k}|}{d}.\qquad\qquad\qquad\qquad\;\nonumber
\end{eqnarray}
Combining the above analysis, the inequality holds
\begin{eqnarray}
|\mathbb{P}(M_i^k\in\mathcal{O}^k|y_i^k)-\mathbb{P}(M_i^{'k}\in\mathcal{O}^k|y_i^{'k})|\nonumber\\
\leq\left(\frac{\sigma}{1-a}\sqrt{\frac{2c}{e\gamma^k}}+\frac{1}{d}\right)|y_i^k-y_i^{'k}|.\label{case1}
\end{eqnarray}

\noindent\textbf{Case II:  \(\bm{\mathbf{\varsigma_i^k=\varsigma_i^{'k}=0}}\).} In this case, player $i$ does not send any message to its neighbors at iteration $k$, i.e., $M_i^k=\emptyset$. Then, we calculate
\begin{eqnarray}
|\mathbb{P}(M_i^k=\emptyset)-\mathbb{P}(M_i^{'k}=\emptyset)|\qquad\qquad\qquad\qquad\nonumber\\
=|1-\mathbb{P}(\varsigma_i^k=1|y_i^k)-(1-\mathbb{P}(\varsigma_i^{'k}=1|y_i^{'k}))|\nonumber\\
=|\mathbb{P}(\varsigma_i^k=1|y_i^k)-\mathbb{P}(\varsigma_i^{'k}=1|y_i^{'k})|\qquad\qquad\; \nonumber\\
\leq \frac{\sigma}{1-a}\sqrt{\frac{2c}{e\gamma^k}}|y_i^k-y_i^{'k}|.\qquad\qquad\qquad\qquad \label{case2}
\end{eqnarray}
Summarizing the results in Case I and Case II, one has
\begin{eqnarray}
|\mathbb{P}(M_i^k\in\mathcal{O}^k)-\mathbb{P}(M_i^k\in\mathcal{O}^k)|\qquad \nonumber\\
\leq\left(\frac{\sigma}{1-a}\sqrt{\frac{2c}{e\gamma^k}}+\frac{1}{d}\right)|y_i^k-y_i^{'k}|.
\end{eqnarray}
Based on the composition theorem for differential privacy in \cite{Dwork}, it follows from the conditions of Theorem 2 that
\begin{eqnarray}
\sup_{\mathrm{Adj}(\mathbf{F},\mathbf{F^{'}})}|\mathbb{P}(M^k(\mathbf{F})\in \mathcal{O}^k)-\mathbb{P}(M^k(\mathbf{F^{'}})\in\mathcal{O}^k)| \nonumber\\
\leq \sum_{i=1}^N\left(\frac{\sigma}{1-a}\sqrt{\frac{2c}{e\gamma^k}}+\frac{1}{d}\right)|y_i^k-y_i^{'k}|\,\nonumber\\
\leq \left(\frac{\sigma}{1-a}\sqrt{\frac{2c}{e\gamma^k}}+\frac{1}{d}\right)\Delta^k\qquad\qquad\ \,\nonumber\\
\leq\left(\frac{\sigma}{1-a}\sqrt{\frac{2c}{e\gamma^k}}+\frac{1}{d}\right)\tilde{C}\frac{(\lambda^k)^2}{\gamma^k}=\delta^k,
\end{eqnarray}
which means that Algorithm 1 is $(0,\delta^k)$-differentially private at each iteration $k$.

With the help of \cite{Dwork}, if Algorithm 1 is $(0,\delta^k)$-differential privacy at each iteration $k$, then this algorithm for the composition of $T$ iteration can satisfy $(0,\sum_{k=0}^T\delta^k)$-differential privacy. It can be directly inferred that 
\begin{eqnarray}
\sum_{k=0}^T\delta^k=\sum_{k=0}^T\left(\frac{\sigma}{1-a}\sqrt{\frac{2c}{e\gamma^k}}+\frac{1}{d}\right)\tilde{C}\frac{(\lambda^k)^2}{\gamma^k}.
\end{eqnarray}
Obviously, if the sequence $\{\frac{(\lambda^k)^2}{(\gamma^k)^{3/2}}\}$ is summable, then we can ensure that $\sum_{k=0}^T\delta^k$ is finite for $T\in\infty$ and $\sum_{k=0}^T\delta^k<1$ can be satisfied by selecting sufficiently large quantization interval and sufficiently small trigger threshold tuning coefficient.


 \noindent\textbf{Remark 6.} As a pioneering work, a ternary quantization-based privacy-preserving mechanism was proposed in \cite{Wang2022} for  distributed optimization, which requires the design parameter of the quantizer not exceed the $l_{\infty}$ norm of its input. This requirement, however, renders the mechanism inapplicable to the distributed NE seeking problem, where the quantizer's input is the decision estimate whose $l_{\infty}$ norm cannot be predetermined. In contrast, the parameters of our  privacy-preserving mechanism (\ref{trigger})-(\ref{quantizer}) are independent of players' decisions and their estimates, and can be set in advance. 

\section{Simulation}

 In this section, an energy consumption game  is considered to demonstrate the effectiveness of Algorithm 1, where the communication graph for players is shown in Fig. 1. The players' cost functions of the game are $f_i(\mathbf{x})=(x_i-\hat{x}_i)^2+(0.04\sum_{i=1}^5x_i+5)x_i$ where $\hat{x}_1=50$, $\hat{x}_2=55$, $\hat{x}_3=60$, $\hat{x}_4=65$, and $\hat{x}_5=70$. Noting that each player $i$ has limited capability to adjust its energy consumption, the decision variables are constrained as $x_1\in[40,45]$, $x_2\in[44,49]$, $x_3\in[48,53]$, $x_4\in[54,59]$, and $x_5\in[58,63]$. 
It follows from $\nabla_if_i(\mathbf{x})=0$ that the unique NE is given as $\mathbf{x}^*=[41.5,46.4,51.3,56.2,61.1]^T$.
\begin{figure}
    \begin{center}
    \includegraphics[height=3cm]{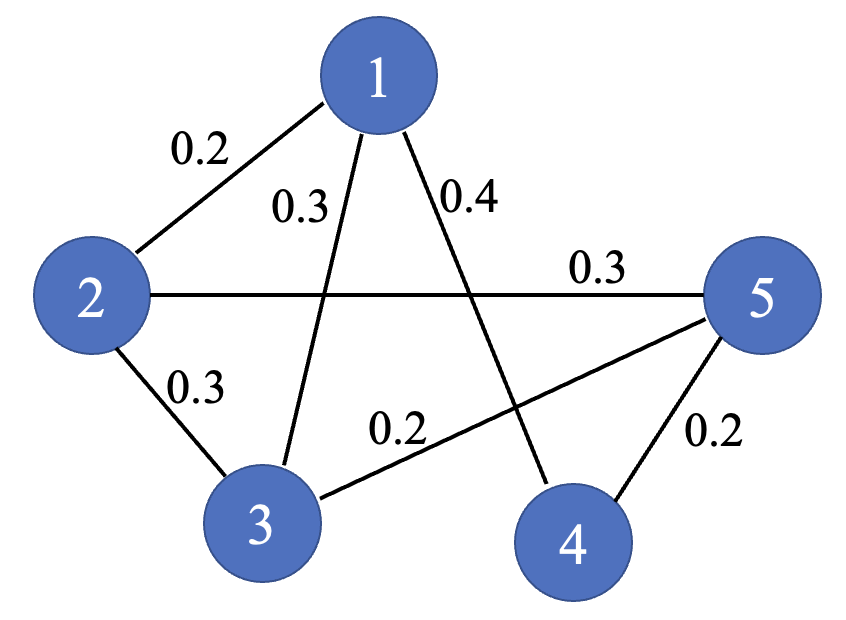}    
    \caption{Communication graph for players.} 
    \label{fig1}                                 
    \end{center}                                 
 \end{figure}
\begin{figure}
    \begin{center}
    \includegraphics[height=4.8cm]{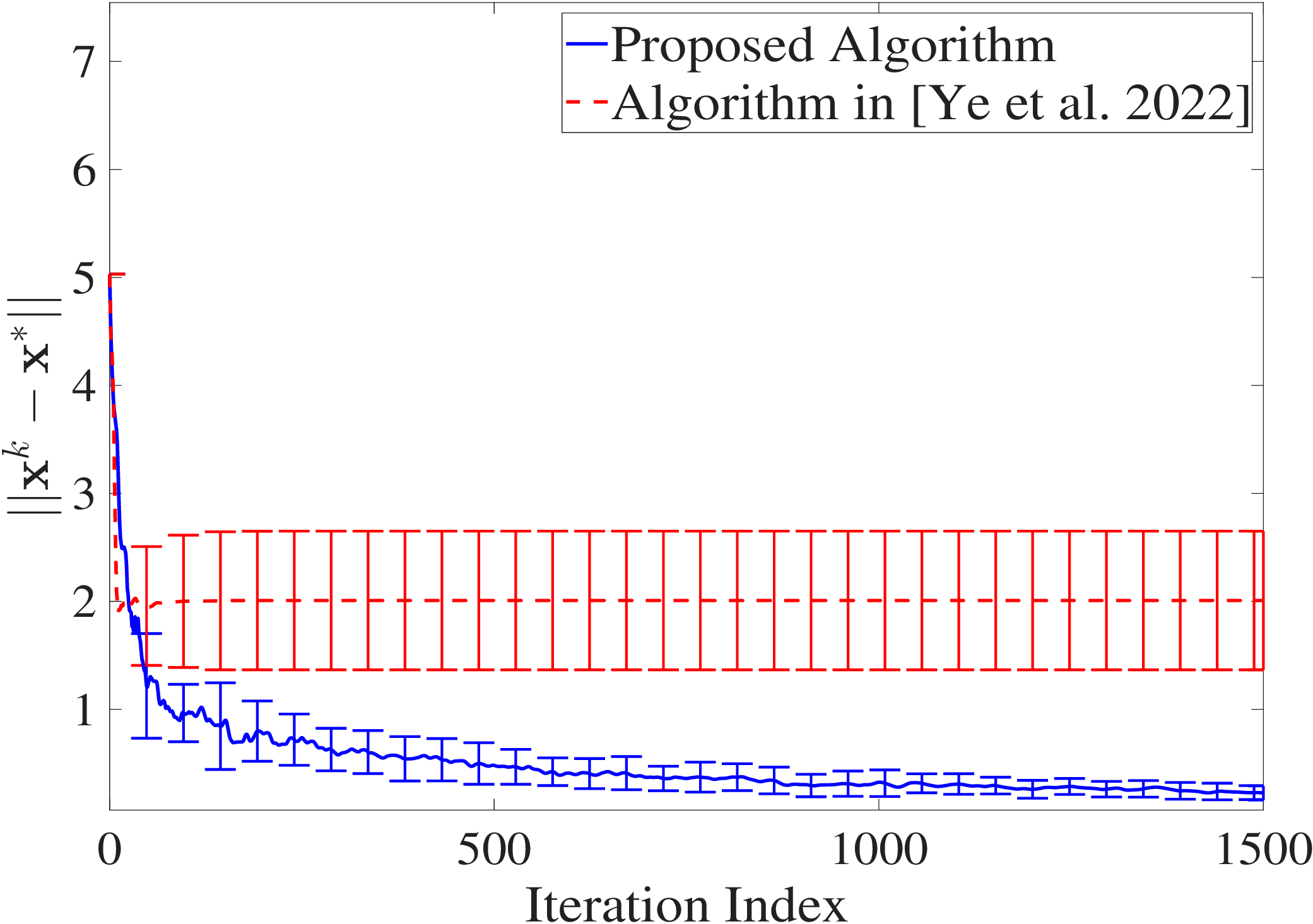}    
    \caption{Comparison of Algorithm 1 with the differential-privacy approach for distributed aggregative games by \cite{Ye2021}.} 
    \label{fig1}                                 
    \end{center}                                 
 \end{figure}
  \begin{figure}
    \begin{center}
    \includegraphics[height=4.7cm]{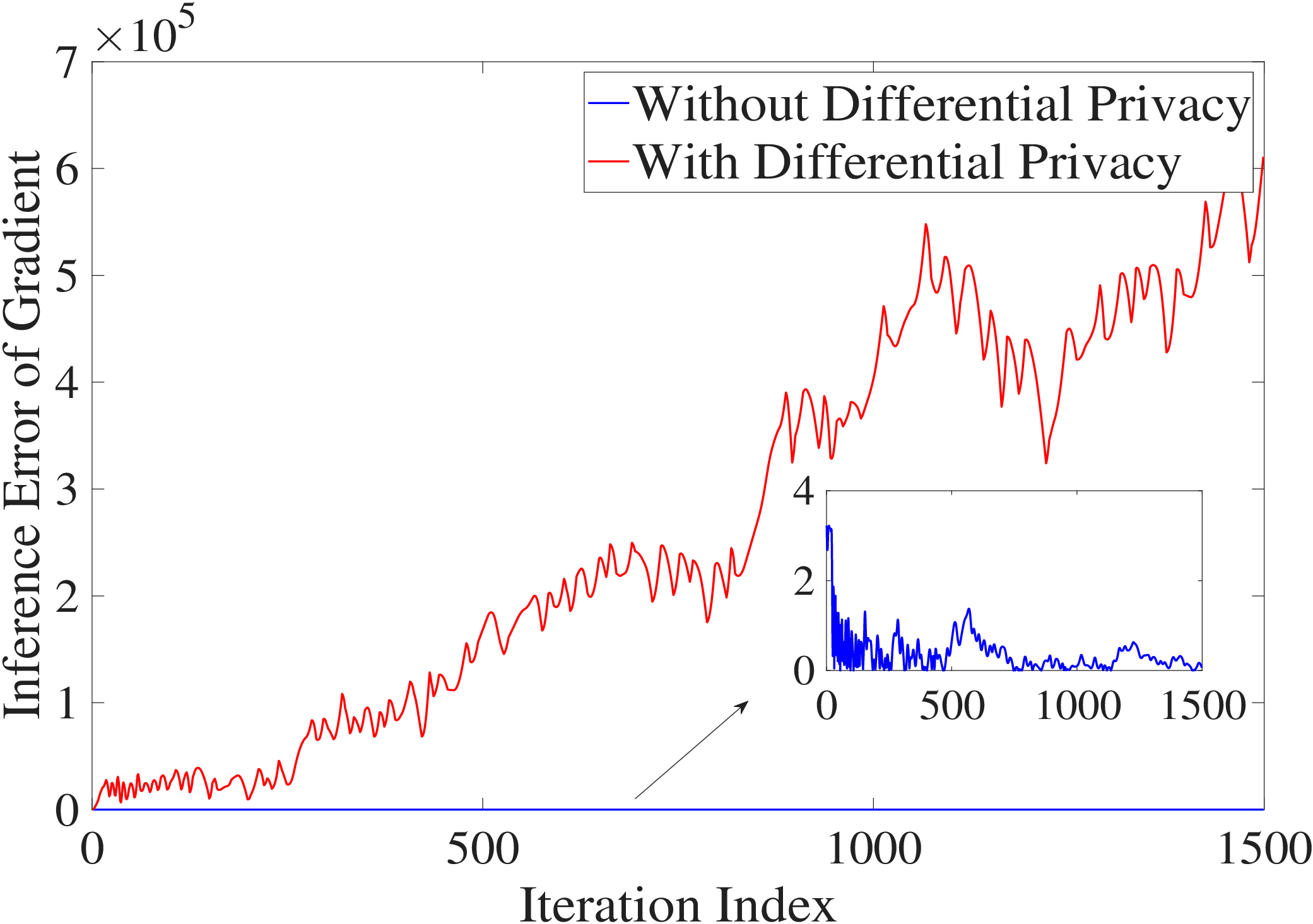}    
    \caption{Inference error of gradient with and without differential privacy.} 
    \label{fig1}                                 
    \end{center}                                 
 \end{figure}
\begin{figure}
    \begin{center}
    \includegraphics[height=4.8cm]{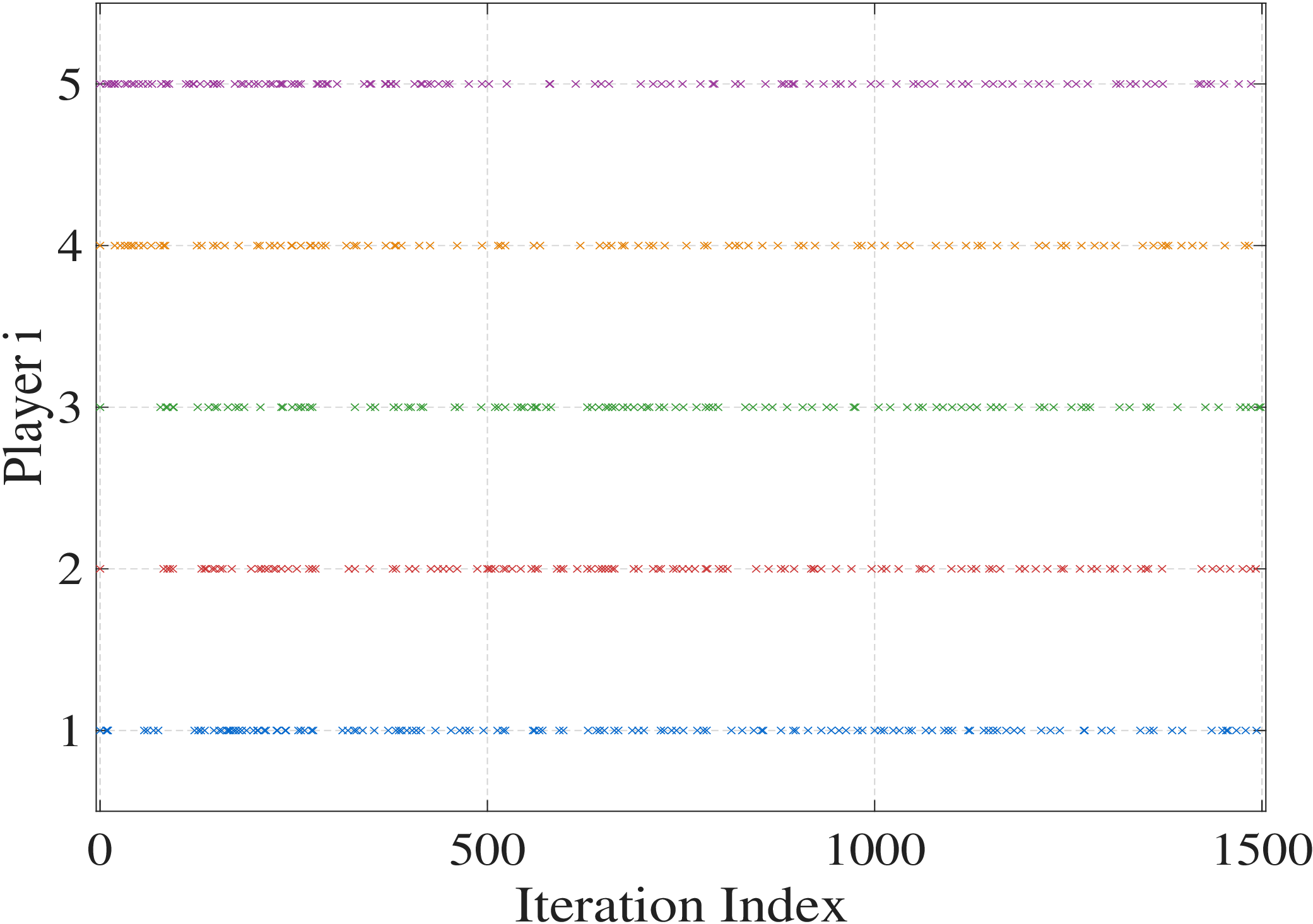}    
    \caption{Triggering times for each player.} 
    \label{fig1}                                 
    \end{center}                                 
 \end{figure}
To protect the sensitive information of players in the game, we use the privacy-preserving mechanism in the form of  (\ref{trigger})-(\ref{quantizer}), in which the quantization interval is chosen as $d=15$ and the triggering parameters are set as $\sigma=1.03$, $a=0.05$, and $c=0.0001$. To evaluate the performance of the proposed Algorithm 1, the stepsize $\lambda^k$ and decaying factor sequence $\gamma^k$ are designed as $\lambda^k=\frac{0.03}{1+0.01k^{0.95}}$ and $\gamma^k=\frac{1.2}{1+0.12k^{0.55}}$, which satisfy the conditions in Theorems 1 and 2. In simulation, the algorithm is run for $1500$ times and it is not difficult to deduce that $\delta^{1500}=0.046$, i.e., Algorithm 1 is $(0,0.046)$-differential privacy at iteration $k=1500$. The average (shown by the curve) as well as the variance (shown by the error bars) of the distance $\left\|\mathbf{x}^k-\mathbf{x}^*\right\|$ between the decision variable $\mathbf{x}^k$ and the NE $\mathbf{x}^*$ are shown in Fig. 2. For comparison, we also run the existing differential-privacy distributed NE seeking algorithm proposed in \cite{Ye2021}. It is worth noting that a geometrically decreasing stepsize was used in the differential privacy method of \cite{Ye2021} to protect privacy,  but this rapidly decreasing stepsize leads to an imprecise convergence to the NE. It can be clearly seen from Fig. 2 that the distance $\left\|\mathbf{x}^k-\mathbf{x}^*\right\|$ of the algorithm in \cite{Ye2021} will no longer decrease after running approximately $200$ times, while our algorithm has a comparable convergence speed but much better accuracy.  Noting that the sensitive information of each player is their cost function, we evaluate the inference error of the gradient of player 1's cost function to verify the effectiveness of the proposed differential privacy algorithm. Specifically, we compute the absolute error between the inferred gradient and the true gradient under Algorithm 1, both with and without differential privacy. As shown in Fig. 3, our algorithm with differential privacy yields a significantly larger inference error. This means that it is very difficult for eavesdroppers to infer the original private information from the observed data sequence,  which demonstrates the effectiveness of our differentially private algorithm. Fig. 4 shows the triggering times for each player and the corresponding trigger probabilities are $9.19\%$, $8.33\%$, $7.53\%$, $6.93\%$, and $8.79\%$, respectively.  This indicates that continuous communication is successfully avoided under (\ref{trigger}), which is beneficial both for privacy protection and communication efficiency.

\section{Conclusion}
This article studies the privacy-preserving distributed NE seeking for aggregative games. A novel differential privacy mechanism is developed, which achieves strong privacy protection by randomizing both the timing and the content of the information exchange. A distributed NE seeking algorithm with dual randomness is proposed. We theoretically prove that the algorithm achieves rigorous $(0,\delta)$-differential privacy at each iteration while preserving provable convergence accuracy. More importantly, in contrast to existing works that fail to maintain privacy over long-term operation, our algorithm sustains this privacy guarantee over an infinite number of iterations, while saving communication resource. 

\bibliographystyle{plain}        

\end{document}